\documentstyle[12pt,amscd]{amsart}
\newtheorem{Theorem}{Theorem}[section] \newtheorem{Lemma}[Theorem]{Lemma}
\newtheorem{Corollary}[Theorem]{Corollary} \newtheorem{Proposition}[Theorem]{Proposition}
\newtheorem*{Definition}{Definition}

\def\Ext{\operatorname{Ext}}\def\Tor{\operatorname{Tor}}\def\grade{\operatorname{grade}}\def\Coker{\operatorname{Coker}}\def\depth{\operatorname{depth}}\def\height{\operatorname{ht}}\def\proj{\operatorname{proj.dim}}\def\Ann{\operatorname{Ann}}\def\rank{\operatorname{rank}}\def\id{\operatorname{id}}\let\To=\longrightarrow\def\Sing{\operatorname{Sing}} \def\NCM{\operatorname{N_{CM}}}  
\def\sk{\par\smallskip} \def\md{\par\medskip}

\begin{document}
\title{specialization of modules over a local ring}
\author{Dam Van Nhi\and Ng\^o Vi\^et Trung}
\address{Pedagogical College (CDSP), Thai Binh, Vietnam} \address{Institute of Mathematics, Box 631, B\`o H\^o, Hanoi, Vietnam} \thanks{The authors are partially  supported by the National Basic Research Program}\maketitle 
\centerline{\em Dedicated to David Buchsbaum}\bigskip

\section*{Introduction}\md

Given an object defined for a family of parameters $u=(u_1,\ldots,u_m)$ we can often substitute $u$ by a family $\alpha=(\alpha_1,\ldots,\alpha_m)$ of elements of an infinite field $K$ to obtain a similar object  which is called a specialization. The new object usually behaves like the given object  for almost all $\alpha$, that is, for all $\alpha$ except perhaps those lying on a proper algebraic subvariety of $K^m.$  Though specialization is a classical method in Algebraic Geometry, there is no systematic theory for what can be ``specialized''.\sk

The first step toward an algebraic theory of specialization was the introduction of the specialization of an ideal by W. Krull [K1],  [K2]. Given an ideal $I$ in a polynomial ring $R=k(u)[X]$, where $k$ is a subfield of $K$, he defined the specialization of $I$ as the ideal $$I_\alpha=\{f(\alpha,X)|\ f(u,X)\in I\cap k[u,X]\}$$ of the polynomial ring $R_\alpha=k(\alpha)[X]$. For almost all $\alpha\in K^m$, $I_{\alpha}$ inherits most of the basic properties of $I$. A. Seidenberg [S] and W.-E. Kuan [Ku1], [Ku2], [Ku3] used specializations of ideals to prove that hyperplane sections of normal varieties are normal again under certain conditions. Inspired by their works, the second author studied  the preservation of  properties from the factor ring $R/I$ to $R_\alpha/I_\alpha$ with regard to hypersurface sections [T2]. However, since there was no appropriate theory for the specialization of modules, many ring-theoretic properties that are characterized by homological means could not be investigated.\sk

The task of finding a suitable definition for specializations of modules is not so easy since we can not specialize elements of an arbitrary $R$-modules. Recently, using finite free presentations the authors introduced the specialization of an $R$-module which covers Krull's specialization of an ideal and for which a systematic theory can be developed [NT]. The aim of this paper is to define and to study specializations of modules over a local ring $(R/I)_P$, where $P$ is an arbitrary prime ideal of $R$. For this we only need to consider $R_P$-modules, and we will follow the approach of the previous paper. There are two obstacles for the task. Firstly, the specialization $P_\alpha$ of $P$ needs not to be a prime ideal so that we may have different specializations for a single local ring $R_P$. Secondly, unlike the case of $R$-modules, there is no existing theory on specializations of ideals of $R_P$. Despite these obstacles we can develop a systematic theory of specializations of $R_P$-modules and we can show that almost all properties of $R_P$-modules specialize.\sk

The paper is divided in four sections. In Section 1 we introduce specializations of finitely generated $R_P$-modules and specializations of homomorphisms between them. We shall see that these notions are unique up to (canonical) isomorphisms. Moreover, we also show that the new notions are compatible with the specializations of $R$-modules introduced in the previous paper. In Section 2 we study basic properties of specializations of modules. The key result is the preservation of exact sequences. This result allows us to show that specializations preserve submodules and ideals and basic operations between them. In particular, specializations do not change the dimension, the height and, if $k$ is a perfect field, also the length. In Section 3 we prove that specializations commute with the Tor and Ext functors. Using Grothendieck's local duality we can show that the annihilators of local cohomology modules specialize. As consequences, specializations preserve the Cohen-Macaulayness, the generalized Cohen-Macaulayness, and, if $k$ is a perfect field, also the Buchsbaumness. Section 4 is devoted to specializations of factor rings of $R_P$. We will study when such a specialization is a domain and  we shall see that it preserves the Gorensteiness and the normality. \sk

We will keep the notations of this introduction throughout this paper and refer to [BH] and [E] for unexplained terminology. For convenience we will often skip the phrase ``for almost all $\alpha$'' when we are working with specializations.\md

\section{Definition of specializations of $R_P$-modules}\md

Let $P$ be an arbitrary prime ideal of $R$. The first obstacle in defining the specialization of $R_P$ is that the specialization $P_\alpha$ of $P$ need not to be a prime ideal.\md

\noindent {\em Example.} Let $P=(uX^2 -1)\subseteq {\Bbb C}(u)[X].$  Then $P$ is a prime ideal, whereas $$P_\alpha=(\alpha X^2-1)=(\sqrt{\alpha}X-1)\cap(\sqrt{\alpha}X+1)$$ is not a prime ideal  for all $\alpha\in\Bbb C$.\md

The natural candidate for the specialization of $R_P$ is the local ring $(R_\alpha)_\wp,$ where $\wp$ is an arbitrary associated prime ideal of $P_\alpha.$ Such a local ring was already considered with regard to specializations of points in [T2].\sk

\begin{Definition}  {\rm We call $(R_\alpha)_\wp$ a {\it specialization} of $R_P$ with respect to $\alpha.$}\end{Definition}\sk

For short we will put $S=R_P$ and $S_\alpha=(R_\alpha)_\wp.$  The notion $S_\alpha$ is not unique. However, all local rings $S_\alpha$ have the same dimension as $S.$ This follows from the case $I=P$ of the following more general fact. 

\begin{Lemma} {\rm [T2, Lemma 1.5]} Let $I$ be an ideal in $R.$ For almost all $\alpha,$\sk
{\rm (i)} $I_\alpha$ is unmixed if $I$ is unmixed,\sk
{\rm (ii)} $\height I_\alpha S_\alpha=\height IS.$\end{Lemma}

Let $f$ be an arbitrary element of $R.$ We may write $f= p(u,x)/q(u),\; p(u, x)\in k[u, x],\ q(u)\in k[u]\setminus\{0\}.$  For any $\alpha$ such that $q(\alpha )\ne 0$ we define $f_\alpha:= p(\alpha,x)/q(\alpha).$ It is easy to check that this element does not depend on the choice of $p(u,x)$ and $q(u)$ for almost all $\alpha.$ Now, for every fraction $a=f/g,\; f,g\in R,\ g\ne 0,$  we define $a_\alpha:=f_\alpha/g_\alpha$ if $g_\alpha\ne 0.$ Then $a_\alpha$ is uniquely determined for almost all $\alpha.$\sk

The following lemma shows that the above definition of $S_\alpha$ reflects the intrinsic substitution $u\to\alpha$ of elements of $R.$

\begin{Lemma} Let $a$ be an arbitrary element of $S.$ Then $a_\alpha\in S_\alpha$ for almost all $\alpha.$\end{Lemma}

\begin{pf} Let $a=f/g$ with $f,g\in R$, $g\notin P$. We have to show that $g_\alpha\notin\wp.$ Since $P$ is prime, $P: g=P.$ By [Kr1, Satz 3], $P_\alpha=(P:g)_\alpha=P_\alpha: g_\alpha$. Hence $g_\alpha\notin\wp$. \end{pf}
 
To define the specialization of a finitely generated $S$-module we need to define the specialization of a homomorphism of free $S$-modules of finite ranks.\sk 

Let $F, G$ be finitely generated free $S$-modules. Let $\phi: F\To G $ be an arbitrary  homomorphism of free $S$-modules of finite ranks. With fixed bases of $F$ and $G$,  $\phi$ is given by a  matrix $A=(a_{ij}),\ a_{ij} \in S.$ By Lemma 1.2, the matrix $A_\alpha :=((a_{ij})_\alpha)$ has all its entries in $S_\alpha$ for almost all $\alpha.$ Let $F_\alpha$ and $G_\alpha$ be free $S_\alpha$-modules of the same rank as $F$ and $G,$ respectively.

\begin{Definition} {\rm For fixed bases of $F_\alpha$ and $G_\alpha,$ the homomorphism $\phi_\alpha: F_\alpha\To G_\alpha$ given by the matrix $A_\alpha$ is called the {\it specialization} of $\phi$ with respect to $\alpha.$}\end{Definition}

The definition of $\phi_\alpha$ does not depend on the choice of the bases of $F, G$ in the sense that if $B$ is the matrix of $\phi$ with respect to other bases of $F,G,$ then there are bases of $F_\alpha,G_\alpha$ such that $B_\alpha$ is the matrix of $\phi_\alpha$ with respect to these bases. 

\begin{Lemma} Let $\phi,\psi: F\To G$ and $\delta: G\To H$ be homomorphisms of finitely generated free $S$-modules. Then, for almost all $\alpha$, \begin{eqnarray*}  (\phi +\psi)_\alpha &=&\phi_\alpha +\psi_\alpha,\\ (\delta\phi)_\alpha &=&\delta_\alpha\phi_\alpha.\end{eqnarray*}\end{Lemma}

\begin{pf} This follows from the matrix presentations of the homomorphisms.\end{pf}

Now we define the specialization of an arbitrary finitely generated $S$-module as follows. 

\begin{Definition} {\rm Let $L$ be a finitely generated $S$-module and $F_1\overset\phi\To F_0\To L\To 0$ a finite free presentation of $L.$ The $S_\alpha$-module $L_\alpha:=\Coker\phi_\alpha$ is called a {\it specialization} of $L$ (with respect to $\phi$).}\end{Definition}

\noindent{\em Remark.}  If $L$ is a free $S$-module of finite rank, then $0\To L\To L\To 0$ is a finite free presentation of $L.$ Therefore $L_\alpha$ as defined above is a free $S_\alpha$-module with $\rank L_\alpha=\rank L.$ This shows that the earlier introduction of $F_\alpha$ and $G_\alpha$ is compatible with the above definition.\md

The definition of $L_\alpha$ depends on the chosen presentation of $L.$ To show that $L_\alpha$ is uniquely determined up to isomorphisms we need to define the specialization of a homomorphism of finitely generated $S$-modules.\sk

Let $v: L\To M$ be a homomorphism of finitely generated $S$-modules. Consider a commutative diagram $$\CD F_1 @>\phi>> F_0 @>>>L @>>> 0\\ @VVv_1 V   @VVv_0 V  @VVv V\\ G_1 @ >\psi>>G_0 @>>>M @>>> 0,\\\endCD$$ where the rows are finite free presentations of $L,M.$ By Lemma 1.3, $\psi_\alpha (v_1)_\alpha=(v_0)_\alpha\phi_\alpha.$ Hence there is an induced homomorphism $v_\alpha: L_\alpha\To M_\alpha,$ which  makes the diagram 
$$\CD (F_1)_\alpha @>\phi_\alpha >>(F_0)_\alpha @>>>L_\alpha @>>> 0\\ @VV(v_1)_\alpha V @VV(v_0)_\alpha V @VVv_\alpha V\\ (G_1)_\alpha @ >\psi_\alpha >>(G_0)_\alpha @>>>M_\alpha @>>> 0\endCD$$ commutative for almost all $\alpha.$\sk

\begin{Definition}  {\rm We call the induced homomorphism $v_\alpha$ a {\it specialization} of $v: L\To M$ with respect to $(\phi,\psi).$}\end{Definition}

The above definition does not depend on the choice of $v_0,v_1.$ Indeed, if we are given two other maps $w_i: F_i\To G_i,\;i=0,1,$ which lift the same homomorphism $v: L\To M,$ then the maps $(w_i)_\alpha$ induce the same map $v_\alpha: L_\alpha\To M_\alpha.$ In particular, the specialization of $\id_L$ with respect to $(\phi,\phi)$ is the identity map $\id_{L_\alpha}.$

\begin{Lemma} Let $v,w: L\To M$ and $u: M\To N$ be homomorphisms of finitely generated $S$-modules. Then, for almost all $\alpha$, \begin{eqnarray*} (v +w)_\alpha &=& v_\alpha +w_\alpha,\\  (u v)_\alpha &=& u_\alpha v_\alpha.\end{eqnarray*}\end{Lemma}

\begin{pf} This is an easy consequence of Lemma 1.2 and the above definition.\end{pf}

Now we will see that $L_\alpha$ is uniquely determined up to canonical isomorphisms. Let $F_1\overset\phi\To F_0\To L\To 0$ and $G_1\overset\psi\To G_0\To L\To 0$ be two different presentations of $L.$ Let $\varepsilon_1: L_\alpha\To L'_\alpha$ and  $\varepsilon_2: L'_\alpha\To L_\alpha$ denote the specialization of $\id_L$ with respect to $(\phi,\psi)$ and to $(\psi,\phi),$ respectively. By Lemma 1.4, $\varepsilon_2\varepsilon_1$ and $\varepsilon_1\varepsilon_2$ are the specializations of  $\id_L$ with respect to $(\phi,\phi)$ and $(\psi,\psi).$ Hence  $\varepsilon_2\varepsilon_1=\id_{L_\alpha}$ and $\varepsilon_1\varepsilon_2=\id_{L'_\alpha}.$ Therefore $\varepsilon_1=(\id_L)_\alpha$ is an isomorphim.\par

We end this section by showing that our definition of specializations of $S$-modules is consistent with that of $R$-modules introduced in [NT].

\begin{Proposition} Let $D$ be a finitely generated $R$-module. For almost all $\alpha$, $$(D\otimes_R S)_\alpha\cong D_\alpha\otimes_{R_\alpha} S_\alpha.$$\end{Proposition}
 
\begin{pf} Let  $E_1\overset {\varphi}\To E_0\To D\To 0$ be a finite free presentation of $D.$ Let $A=(a_{ij})$ be a representing matrix for $\varphi.$ Tensoring with $S$ we get a finite free presentation: $$E_1\otimes_R S\overset {\varphi\otimes S}\To E_0\otimes_R S\To D\otimes_R S\To 0,$$ where $\varphi\otimes S$ is represented by the matrix $A.$ By the definition of specialization of $S$-module, $(D\otimes S)_\alpha$ is the cokernel of the map $(\varphi\otimes S)_\alpha:(E_1\otimes_R S)_\alpha\To (E_0\otimes_R S)_\alpha,$ where $(\varphi\otimes S)_\alpha$ is represented by the matrix $A_\alpha.$ On the other hand, from the exact sequence $(E_1)_\alpha\overset {\varphi_\alpha}\To (E_0)_\alpha\To D_\alpha\To 0$ we get the exact sequence $$(E_1)_\alpha\otimes_{R_\alpha} S_\alpha\overset {\varphi_\alpha\otimes S_\alpha}\To (E_0)_\alpha\otimes_{R_\alpha} S_\alpha\To D_\alpha\otimes_{R_\alpha} S_\alpha\To 0,$$ where $(E_1)_\alpha\otimes_{R_\alpha} S_\alpha$ and $(E_0)_\alpha\otimes_{R_\alpha} S_\alpha$ are free $S_\alpha$-modules of the same ranks as $(E_1\otimes_R S)_\alpha$ and $(E_0\otimes_R S)_\alpha,$ respectively, and  $\varphi_\alpha\otimes S_\alpha$ is represented by the matrix $A_\alpha.$ Thus, we can conclude that $(D\otimes_R S)_\alpha\cong D_\alpha\otimes_{R_\alpha} S_\alpha.$\end{pf}

\begin{Corollary} Let $I$ be an ideal of $R.$ Then $(IS)_\alpha\cong I_\alpha S_\alpha$ for almost all $\alpha.$\end{Corollary}

\begin{pf} Note that $I_\alpha$ is also the specialization of $I$ as an $R$-module [NT]. Since $IS=I\otimes_R S$ and $I_\alpha S_\alpha=I_\alpha\otimes_{R_\alpha} S_\alpha,$ the conclusion follows from Proposition 1.5.\end{pf}\sk 

\section{Basic properties of specializations of $R_P$-modules}\md 

We will first show that specializations of  $S$-modules ($S = R_P$) preserve the exactness of finite complexes of free modules. \sk

Let $\phi: F\To G$ is a homomorphism of free $S$-modules of finite ranks and let $A$ be the matrix of $\phi$ with respect to fixed bases of $F,G$. Denote by $I_t (\phi)$ the ideal generated by the $t\times t$-minors of $A$. Define $\rank\phi=\max\{t|\, I_t (\phi)\ne 0\}$ and put $I(\phi):=I_d (\phi)$ if $d=\rank\phi.$\sk

Let ${\bold F_\bullet}:\; 0\To F_\ell\overset{\phi_\ell}\To  F_{\ell-1}\To\cdots\To F_1\overset{\phi_1}\To F_0 $ be a complex of free $S$-modules of finite ranks. By a criterion of Buchsbaum and Eisenbud, ${\bold F_\bullet}$ is exact if and only if
\begin{eqnarray*}\rank\phi_j +\rank\phi_{j +1} &=&\rank F_j,\\\depth I(\phi_j) &\ge& j\end{eqnarray*}
for $j=1,\ldots,\ell$ (see [BE, Corollary 1]). By Lemma 1.3 the sequence $$({\bold F_\bullet})_\alpha:\; 0\To(F_\ell)_\alpha\overset {(\phi_\ell)_\alpha}\To (F_{\ell-1})_\alpha\To\cdots\To (F_1)_\alpha\overset{(\phi_1)_\alpha}\To (F_0)_\alpha$$ is a complex of free $S_\alpha$-modules. Hence we can also apply the above criterion to $({\bold F_\bullet})_\alpha$.

\begin{Lemma} Let ${\bold F_\bullet}$ be an exact complex of free $S$-modules of finite ranks. Then $({\bold F_\bullet})_\alpha$ is an exact complex for almost all $\alpha.$\end{Lemma}

\begin{pf} By definition, $\rank (F_j)_\alpha=\rank F_j$ for $j=0,\ldots,\ell.$ It is easy to check that $\rank (\phi_j)_\alpha=\rank\phi_j$ (see [NT, Lemma 1.4 (i)]). Hence $$\rank (\phi_j)_\alpha +\rank (\phi_{j+1})_\alpha = \rank\phi_j +\rank\phi_{j+1}= \rank F_j\;=\;\rank (F_j)_\alpha.$$ It remains to show that $\depth I((\phi_j)_\alpha)\ge j.$ Let $\phi_j$ be given by a matrix $A_j.$ Let $a_1,\ldots,a_s$ be the $d\times d$-minors of $A_j,d=\rank\phi_j.$ We may write $a_j=f_j/g_j$ with $f_j,g_j\in R$, $g_j\not\in P.$ Then  $I((\phi_j)_\alpha)=((a_1)_\alpha,\ldots,(a_s)_\alpha)=((f_1)_\alpha,\ldots,(f_s)_\alpha)S_\alpha.$ Let $I$ be the ideal of $R$ generated by the elements $f_1,\ldots,f_s.$ Then   $I_\alpha=((f_1)_\alpha,\ldots,(f_s)_\alpha)$ by  [NT, Corollary 3.3]. Therefore, $I((\phi_j)_\alpha)=I_\alpha S_\alpha$. By Lemma 1.1, $\height I_\alpha S_\alpha=\height IS.$ Note that $IS = I(\phi_j)$. Then, since $S$ and $S_\alpha$ are Cohen-Macaulay rings,\md
$\quad\quad\quad \depth I((\phi_j)_\alpha) = \height I((\phi_j)_\alpha) = \height I_\alpha S_\alpha = \height I(\phi_j) = \depth  I(\phi_j)\ge j.$ \end{pf}  

Now we can proceed as in [NT] to establish basic properties of specializations of $S$-modules.\sk

\begin{Theorem} Let $\;0\To L\To M\To N\To 0$ be an exact sequence of finitely generated $S$-modules. Then $0\To L_\alpha\To M_\alpha\To N_\alpha\To 0$ is exact for almost all $\alpha.$\end{Theorem}
 
\begin{pf} See the proof  of [NT, Theorem 2.4]. \end{pf} 

Let $L$ be a finitely generated $S$-module and $M$ an arbitrary submodule of $L$. We can consider $M_\alpha$ as a submodule of $L_\alpha.$ Indeed, the map $M_\alpha\To L_\alpha$ is injective by Theorem 2.2. This map is canonical because different specializations $M_\alpha,\; M'_\alpha$ of $M$ have the same image in $L_\alpha.$ That follows from the commutative diagram $$\begin{array}{rccl} & M_\alpha & &\\ & &\searrow &\\ (\id_M)_\alpha &\wr\Vert & & L_\alpha.\\ & &\nearrow\\ & M'_\alpha & &\end{array}$$ Unless otherwise specified, we will always identify $M_\alpha$ with a submodule of $L_\alpha.$ In particular, the specialization ${\frak a}_\alpha$ of any ideal ${\frak a}$ of $S$ can be identified with an ideal of $S_\alpha.$\sk

\begin{Lemma} Let $M, N$ be submodules of a finitely generated $S$-module $L$. For almost all $\alpha$, there are canonical isomorphisms: \begin{eqnarray*} (L/M)_\alpha &\cong & L_\alpha/M_\alpha,\\ (M\cap N)_\alpha &\cong & M_\alpha\cap N_\alpha,\\ (M+N)_\alpha&\cong& M_\alpha+N_\alpha.\end{eqnarray*}\end{Lemma}

\begin{pf} See the proof of [NT, Theorem 3.2].\end{pf}

\begin{Lemma} Let $L=Se_1 +\cdots +Se_s$ be an $S$-module. For almost all $\alpha,$ there exist elements $(e_1)_\alpha,\cdots,(e_s)_\alpha\in L_\alpha$ such that  $(Se_j)_\alpha=S_\alpha (e_j)_\alpha,\; j=1,\cdots,s,$ and$$L_\alpha=S_\alpha (e_1)_\alpha +\cdots +S_\alpha (e_s)_\alpha.$$
In particular, if $a_1,\ldots,a_s$ are elements of $S$ and ${\frak a} = (a_1,\ldots,a_s)$, then 
$${\frak a}_\alpha = ((a_1)_\alpha,\ldots,(a_s)_\alpha).$$\end{Lemma}

\begin{pf} For the first statement see the proof of [NT, Corollary 3.3]. For the second statement we only need to consider the case ${\frak a} = a_1S$ is a principal ideal. Write $a_1 = f/g$ with $f, g \in R$. Then ${\frak a} = fS$. Note that $(fR)_\alpha = f_\alpha R_\alpha$. Then, using Corollary 1.6 we obtain ${\frak a}_\alpha = (fS)_\alpha = f_\alpha S_\alpha = (a_1)_\alpha S_\alpha.$\end{pf}

\begin{Lemma} Let $L$ be a finitely generated $S$-module and ${\frak a}$ an ideal of $S$. For almost all $\alpha$ we have
\begin{eqnarray*} (0:_L {\frak a})_\alpha&\cong&0:_{L_\alpha}{\frak a}_\alpha,\\ ({\frak a}L)_\alpha&\cong&{\frak a}_\alpha L_\alpha.\end{eqnarray*}\end{Lemma}

\begin{pf} See the proof [NT, Proposition 3.6].\end{pf}

\begin{Theorem} Let $L$ be a finitely generated $S$-module. Then, for almost all $\alpha,$\sk{\rm (i) } $\Ann L_\alpha=(\Ann L)_\alpha,$\sk{\rm (ii) } $\dim L_\alpha=\dim L$.\end{Theorem}

\begin{pf} See the proof [NT, Theorem 3.4]. \end{pf} 

\begin{Corollary} Let ${\frak a}\supseteq\frak b$ be arbitrary ideals of $S.$ For almost all $\alpha$,$$\height {\frak a}_\alpha/{\frak b}_\alpha=\height {\frak a}/{\frak b}.$$\end{Corollary}

\begin{pf}  Let ${\frak a}={\frak q}_1\cap\ldots\cap {\frak q}_s$ be a primary decomposition of $\frak a$. By Lemma 2.3, ${\frak a}_\alpha=({\frak q}_1)_\alpha\cap\ldots\cap ({\frak q}_s)_\alpha$. Hence we only need to consider the case ${\frak a}$ is a primary ideal. Let $\frak p$ denote the associated prime ideal of $\frak a$. Then ${\frak a}_\alpha$ and ${\frak p}_\alpha$ share the same radical. Hence $\height {\frak a}_\alpha/{\frak b}_\alpha = \height {\frak p}_\alpha/{\frak b}_\alpha$. Now we will show that $\height {\frak p}_\alpha/{\frak b}_\alpha =  \height {\frak a}/{\frak b}.$ Let $\wp'$ be an arbitrary associated prime of ${\frak p}_\alpha$. By Proposition 1.6 we may consider $S_{\wp'}$ and ${\frak b}_\alpha S_{\wp'}$ as specialisations of $S_{\frak p}$ and ${\frak b}S_{\frak p}$ with repspect to $\alpha$. By Lemma 2.3,  $(S_{\frak p}/{\frak b}S_{\frak p})_\alpha = S_{\wp'}/{\frak b}_\alpha S_{\wp'}$. Thus, using Theorem 2.6 (ii) we get 
$$\height \wp'_\alpha/{\frak b}_\alpha = \dim S_{\wp'}/{\frak b}_\alpha S_{\wp'} = \dim S_{\frak p}/{\frak b}S_{\frak p} = \height {\frak p}/{\frak b} = \height {\frak a}/{\frak b}.$$  Hence $\height {\frak p}_\alpha/{\frak b}_\alpha =  \height {\frak a}/{\frak b}.$ \end{pf}

\begin{Proposition} Let $L$ be a $S$-module of finite length. Then $L_\alpha$ is a $S_\alpha$-module of finite length for almost all $\alpha$. Moreover, $\ell (L_\alpha)=\ell (L)$ if $k$ is a perfect field. \end{Proposition}

\begin{pf} The first statement follows from the fact that $\dim L_\alpha = \dim L = 0$. To prove the second statement we apply Theorem 2.2 to a composition series of $L$ in order to restrict to the case $\ell(L)=1.$ Then $L\cong S/PS.$ By Corollary 1.6 and Lemma 2.3, $L_\alpha\cong (S/PS)_\alpha=S_\alpha/P_\alpha S_\alpha.$ If $k$ is a perfect field, every extension of $k$ is separable. Thus, $P_\alpha$ is a radical ideal by [Kr2, Satz 14]. So $P_\alpha S_\alpha=\wp S_\alpha$. Hence $\ell(L_\alpha)=\ell (S_\alpha/\wp S_\alpha)=1.$\end{pf} \sk

\section{Preservations of homological properties}\md

This section concerns homological properties which are preserved by specializations. 

\begin{Theorem} Let $L$ be a finitely generated $S$-module. For almost all $\alpha,$\sk
{\rm (i) } $\proj L_\alpha=\proj L,$\sk
{\rm (ii) } $\depth L_\alpha=\depth L.$\end{Theorem}

\begin{pf}  Let ${\bold F_\bullet}:\; 0\To F_\ell\overset{\phi_\ell}\To F_{\ell-1}\To\cdots\To F_1\overset{\phi_1}\To F_0\To L$ be a minimal free resolution of $L.$ Let  $\phi_t$ be given by the matrix $A_t=(a_{tij})$ with respect to fixed bases of $F_0,\ldots,F_\ell,\;t=1,\ldots,\ell.$ Since ${\bold F_\bullet}$ is minimal, $a_{tij}\in PS.$ By Lemma 2.1, the complex 
$$({\bold F_\bullet})_\alpha:\; 0\To (F_\ell)_\alpha\overset{(\phi_\ell)_\alpha}\To(F_{\ell-1})_\alpha\To\cdots\To (F_1)_\alpha\overset{(\phi_1)_\alpha}\To (F_0)_\alpha\To L_\alpha$$ is a free resolution of $L_\alpha.$ Since $(A_t)_\alpha=((a_{tij})_\alpha)$ is the representing matrix of $(\phi_t)_\alpha$ and since $(a_{tij})_\alpha\in P_\alpha S_\alpha=\wp S_\alpha,$ the free resolution $({\bold F_\bullet})_\alpha$ of $L_\alpha$ is minimal. Hence  $\proj L_\alpha=\proj L.$\par
Since $S$ and $S_\alpha$ are Cohen-Macaulay rings, $\depth S_\alpha=\height\wp=\height P=\depth S$ by Lemma 1.1. By Auslander-Buchsbaum formula [E, Theorem 19.9] we have\md
$\quad\quad \depth L_\alpha  = \depth S_\alpha -\proj L_\alpha = \depth S-\proj L\,=\,\depth L.$ \end{pf}
  
\begin{Corollary} Let $L$ be a finitely generated $S$-module. If $L$ is a Cohen-Macaulay $S$-module, then $L_\alpha$ is also a Cohen-Macaulay $S_\alpha$-module for almost all $\alpha.$\end{Corollary}

\begin{pf} By Theorem 2.6 (ii) and Theorem 3.1 (ii) we have\md  $\hspace{4cm} \depth L_\alpha=\depth L=\dim L=\dim L_\alpha.$ \end{pf}

Similarly as in [NT] we can use Theorem 2.2 to show that specializations commute with the Tor and the Ext functors.

\begin{Proposition} Let $L$ and $M$ be finitely generated  $S$-modules. For almost all $\alpha,$ 
\begin{eqnarray*}\Ext^i_{S_\alpha}(L_\alpha,M_\alpha)&\cong&\Ext^i_S(L,M)_\alpha,\;i\ge 0,\\ \Tor^{S_\alpha}_i(L_\alpha,M_\alpha)&\cong&\Tor^S_i(L,M)_\alpha,\; i\ge 0.\end{eqnarray*}\end{Proposition} 

\begin{pf} See the proofs of [NT, Theorem 4.2 and Theorem 4.3].\end{pf}

\begin{Corollary} Let $L$ be a finitely generated $S$-module and $\frak a$ an ideal of $S$. For almost all $\alpha$,\sk
{\rm (i)} $\grade({\frak a}_\alpha,L_\alpha) = \grade({\frak a},L)$,\sk
{\rm (ii)} $\grade L_\alpha = \grade L$.\end{Corollary}

\begin{pf} See the proof of [NT, Corollary 4.4].\end{pf}

Let $\frak m$ and $\frak n$ denote the maximal ideals $PS$ of $S$ and  $\wp S_\alpha$ of  $S_\alpha$. Let $H^i_{\frak m}(L)$ and $H^i_{\frak n}(L_\alpha)$ denote the $i$th local cohomology module of $L$ and $L_\alpha$ with respect to $\frak m$ and $\frak n.$ Put $t=\dim S.$ By the local duality theorem of Grothendieck we have
\begin{eqnarray*} H^i_{\frak m}(L) &\cong&\Ext^{t-i}_S(L,S)^*,\\ \Ext^i_S(L,S)&\cong&H^{t-i}_{\frak m}(L)^*,\end{eqnarray*} where $^*$ denotes the Matlis duality (see e.g. [BH, Section 3.2]). Thus, in order to study $H^i_{\frak m}(L)$ we may pass to $\Ext^{t-i}_S(L,S).$\sk

First we will study the specializations of annihilators of $H^i_{\frak m}(L)$, which carry important information on the struture of $L$ [Fa], [Sc].

\begin{Lemma} Let $L$ be a finitely generated $S$-module. For almost all $\alpha$, $$\Ann H^i_{\frak n}(L_\alpha)=(\Ann H^i_{\frak m}(L))_\alpha,\; i\ge0.$$\end{Lemma}

\begin{pf} By local duality we have \begin{eqnarray*}\Ann H^i_{\frak m}(L)&=&\Ann\Ext^{t-i}_S(L,S)\\ \Ann H^i_{\frak n}(L_\alpha)&=&\Ann\Ext^{t-i}_{S_\alpha}(L_\alpha,S_\alpha).\end{eqnarray*} 
Applying Proposition 3.3 and Theorem 2.6 (i) we obtain $$\Ann\Ext^{t-i}_{S_\alpha}(L_\alpha,S_\alpha)=\Ann\Ext^{t-i}_S(L,S)_\alpha=(\Ann\Ext^{t-i}_S(L,S))_\alpha.$$ Hence $\Ann H^i_{\frak n}(L_\alpha)=(\Ann H^i_{\frak m}(L))_\alpha,\;i\ge 0.$\end{pf}

Recall that a finitely generated module $L$ is said to be a {\it generalized Cohen-Macaulay module} if $\ell (H^i_{\frak m}(L))<\infty$ for all $i=0,\ldots,\dim L-1$ (see e.g. [T3]).

\begin{Theorem} Let $L$ be a  generalized Cohen-Macaulay $S$-module. Then $L_\alpha$ is also a generalized Cohen-Macaulay $S_\alpha$-module for almost all $\alpha.$ Moreover, if $k$ is a perfect field, $$\ell(H^i_{\frak n}(L_\alpha))=\ell(H^i_{\frak m}(L)),\; i=0,\ldots,\dim L-1.$$\end{Theorem}
 
\begin{pf} Let $d=\dim L.$ By Theorem 2.6 (ii), $\dim L_\alpha=d.$  For $i=0,\ldots,d-1,$ $H^i_{\frak m}(L)$ is annihilated by some power of ${\frak m}.$ Hence $H^i_{\frak n}(L_\alpha)$ is annihilated by some power of $\frak n$ by Lemma 3.5. By local duality, $\Ext^{t-i}_{S_\alpha}(L_\alpha,S_\alpha)$  is annihilated by some power of $\frak n,$ too. Since $\Ext^{t-i}_{S_\alpha}(L_\alpha,S_\alpha)$ is a finitely generated $S_\alpha$-module, $\Ext^{t-i}_{S_\alpha}(L_\alpha,S_\alpha)$ has finite length. By [BH, Proposition 3.2.12],  $\ell (H^i_{\frak n}(L_\alpha))=\ell (\Ext^{t-i}_{S_\alpha}(L_\alpha,S_\alpha))<\infty$ for $i=0,\ldots,d-1$. Hence $L_\alpha$ is a generalized Cohen-Macaulay $S_\alpha$-module. 
By Proposition 3.3, $\Ext^{t-i}_{S_\alpha}(L_\alpha,S_\alpha)\cong\Ext^{t-i}_S(L,S)_\alpha.$ If $k$ is a perfect field, we may apply Proposition 2.8 to obtain\md $\hspace{2cm} \ell (H^i_{\frak n}(L_\alpha))=\ell (\Ext^{t-i}_S(L,S)_\alpha)=\ell (\Ext^{t-i}_S(L,S))=\ell (H^i_{\frak m}(L)).$ \end{pf} 

The structure of a generalized Cohen-Macaulay module is best described by standard systems of parameters [T3]. Recall that a system of parameters $a_1,\ldots,a_d$ of $L$ is {\it standard} if $(a_1,\ldots,a_d)H^i_{\frak m}(L/(a_1,\ldots,a_j)L)=0$ for all non-negative integers $i,j$ with $i+j<d$. It may be also characterized by the property that $a^{n_1}_1,\ldots,a^{n_d}_d$ is a $d$-sequence for all positive numbers $n_1,\ldots,n_d$ and for any permutation of $a_1,\ldots,a_d$. An $\frak m$-primary ideal ${\frak a}$ of $S$ is called a {\it standard ideal} for $L$ if every system of parameters of $L$ contained in ${\frak a}$ is standard. In particular, $L$ is a Buchsbaum module if and only if the maximal ideal $\frak m$ is a standard ideal for $L.$ 

\begin{Theorem} Let $L$ be a finitely generated $S$-module and ${\frak a}$ a standard ideal for $L.$ Then ${\frak a}_\alpha$ is a standard ideal for $L_\alpha$ for almost all $\alpha.$\end{Theorem}

\begin{pf} Let $\cal B$ be an $L$-base for ${\frak a},$ that is a generating set for ${\frak a}$ such that every family of $d$ elements of $\cal B$ is a system of parameters of $L$, where $d=\dim L$. It is known that ${\frak a}$ is a standard ideal for $L$ if and only if every family $a_1,\ldots,a_d$ of elements of $\cal B$ is a standard system of parameters of $L$ [T3, Proposition 3.2]. Let $\cal B_\alpha$ denote the set of the specialized elements of $\cal B$ in $S_\alpha.$ By Lemma 2.4, ${\cal B}_\alpha$ is a generating set for ${\frak a}_\alpha.$ By Lemma 2.4 and Lemma 2.5 we have 
$$L_\alpha/((a_1)_\alpha,\ldots,(a_j)_\alpha)L_\alpha=L_\alpha/((a_1,\ldots,a_j)L)_\alpha=(L/(a_1,\ldots,a_j)L)_\alpha$$ for $j=1,\ldots,d$. By Theorem 2.6 (ii), $\dim L_\alpha=d$ and $\dim L_\alpha/((a_1)_\alpha,\ldots,(a_d)_\alpha)L_\alpha=\dim L/(a_1,\ldots,a_d)L=0$. Hence $(a_1)_\alpha,\ldots,(a_d)_\alpha$ is also a system of parameters of $L_\alpha.$ Thus, $\cal B_\alpha$ is an $L_\alpha$-base for ${\frak a}_\alpha.$ By Lemma 3.5,
$$\Ann H_{\frak n}^i(L_\alpha/((a_1)_\alpha,\ldots,(a_j)_\alpha)L_\alpha)=(\Ann H_{\frak m}^i(L/(a_1,\ldots,a_j)L))_\alpha$$ for $i\ge 0$. If $i+j < d$, we have
$(a_1,\ldots,a_d)\subseteq\Ann H_{\frak m}^i(L/(a_1,\ldots,a_j)L).$  Hence
$$((a_1)_\alpha,\ldots,(a_d)_\alpha)H_{\frak n}^i(L_\alpha/((a_1)_\alpha,\ldots,(a_j)_\alpha)L_\alpha)=0.$$ Thus, every family of $d$ elements in ${\cal B}_\alpha$ is a standard system of parameters of $L_\alpha.$ By [T3, Proposition 3.2] we can conclude that ${\frak a}_\alpha$ is a standard ideal for $L_\alpha$.\end{pf} 

\begin{Corollary} Let $L$ be a Buchsbaum $S$-module. Then $L_\alpha$ is a Buchsbaum $S_\alpha$-module for almost all $\alpha$ if  $k$ is a perfect field. \end{Corollary}

\begin{pf} Using Proposition 2.8 we can see that ${\frak m}_\alpha = \frak n$ is the maximal ideal of $S_\alpha$. Since $\frak m$ is a standard ideal for $L$, ${\frak m}_\alpha$ is a standard ideal for $L_\alpha$ by Theorem 3.7. Hence $L_\alpha$ is a Buchsbaum module.\end{pf}\sk

\section{Specialization of factor rings of $R_P$}\md

The aim of this section is to study properties of the factor ring $S_\alpha/{\frak a}_\alpha,$ where $\frak a$ is an arbitrary ideal of $S = R_P.$ For those properties which can be defined for modules we refer to Section 2 and Section 3. Note that $S_\alpha/{\frak a}_\alpha=(S/{\frak a})_\alpha$ by Lemma 2.3.\sk

Following [Kr2] we call an ideal $\frak a$ of $S$ an {\em absolutely prime ideal} if the extension of $\frak a$ in $S\otimes_{k(u)}\overline {k(u)}$ is a prime ideal, where $\overline {k(u)}$ denote the algebraic closure of $k(u).$ 

\begin{Proposition} Let $\frak a$ be an absolutely prime ideal of $S.$ Then ${\frak a}_\alpha$ is a prime ideal for almost all $\alpha.$\end{Proposition}

\begin{pf} Let $I$ be the prime ideal of $R$ such that ${\frak a}=IS.$ Then the extension of $I$ in $\overline {k(u)}[X]$ is also a prime ideal. By [Kr2, Satz 16], $I_\alpha$ is a prime ideal. But ${\frak a}_\alpha=I_\alpha S_\alpha$ by Corollary 1.6. Hence ${\frak a}_\alpha$ is a prime ideal of $S_\alpha.$\end{pf}
  
\begin{Proposition} Let $\frak a$ be an ideal of $S$ such that $S/\frak a$ is a Gorenstein ring. Then $S_\alpha/{\frak a}_\alpha$ is a Gorenstein ring for almost all $\alpha.$\end{Proposition} 

\begin{pf} Let $r=\height\frak a.$ By [BH, Theorem 3.3.7], $S/\frak a$ is a Gorenstein ring if and only if $S/\frak a$ is a Cohen-Macaulay ring and $\Ext^r_S(S/{\frak a},S)\cong S/\frak a$. By Corollary 2.7, $\height{\frak a}_\alpha=r$. By Corollary 3.2, $S_\alpha/{\frak a}_\alpha$ is  a Cohen-Macaulay ring. By Proposition 3.3, 
$$\Ext^r_{S_\alpha}(S_\alpha/{\frak a}_\alpha,S_\alpha)\cong\Ext^r_S(S/{\frak a},S)_\alpha\cong (S/{\frak a})_\alpha\cong S_\alpha/{\frak a}_\alpha.$$Hence $S_\alpha/{\frak a}_\alpha$ is a Gorenstein ring.\end{pf}

Now we will show that specializations preserve the non-Cohen-Macaulay locus and the singular locus. For convenience we denote these locus by $\NCM$ and Sing, respectively. Moreover, for any ideal $I$, we denote by $V(I)$ the set of prime ideals which contain $I$. 

\begin{Lemma} Let $\frak a$ be an arbitrary ideal of $S$. There is an ideal ${\frak b}\supseteq\frak a$ with $\NCM(S/{\frak a})=V({\frak b})$ such that for almost all $\alpha$,$$\NCM(S_\alpha/{\frak a}_\alpha) = V({\frak b}_\alpha).$$ \end{Lemma}

\begin{pf} Let $d = \dim S/\frak a$ and ${\frak a}_i=\Ann H^i_{\frak m}(S/{\frak a})$. Then $\dim S_\alpha/{\frak a}_\alpha= \dim S/{\frak a} = d$ and 
$$({\frak a}_i)_\alpha=(\Ann H^i_{\frak m}(S/{\frak a}))_\alpha=\Ann H^i_{\frak n}((S/{\frak a})_\alpha)=\Ann H^i_{\frak n}(S_\alpha/{\frak a}_\alpha)$$ 
by Lemma 3.5. Put  ${\frak b}=\prod_{i=1}^d\prod_{j > i}({\frak a}_i+{\frak a}_j).$ 
Then $\NCM(S/{\frak a})=V({\frak b})$ by  [Sc, Korollar 6]. By Lemma 2.3 and Lemma 2.5 we have
$${\frak b}_\alpha=\prod_{i=1}^d\prod_{j > i} ({\frak a}_i + {\frak a}_j)_\alpha = \prod_{i=1}^d\prod_{j > i}\big(({\frak a}_i)_\alpha+({\frak a}_j)_\alpha\big).$$ 
Hence $\NCM(S_\alpha/{\frak a}_\alpha)=V({\frak b}_\alpha)$.\end{pf}

\begin{Lemma} Let $k$ be a field of characteristic zero. Let $I$ be an ideal of $R$. There is an ideal $J\supseteq I$ with $\Sing(R/I) = V(J)$ such that for almost all $\alpha$, $$\Sing(R_\alpha/I_\alpha) = V(J_\alpha).$$ \end{Lemma}

\begin{pf} We begin with some observation on the singular locus of $R/I$. Let $r=\height I$ and $A$ a Jacobian matrix of $I$. Let $J_1$ denote the ideal generated by the $r\times r$ minors of $A$. Let $Q\supseteq I$ be a prime ideal such that $(R/I)_Q$ is not regular. If $\height IR_Q=r$, $Q\in V(J_1)$. If $\height IR_Q > r$, $Q$ does not contain any associated prime ideal of $I$ of height $r$. Let $I_r$ denote the intersection of the associated prime ideals of $I$ of height $r$. Let $I'=\cup_{t\ge 1}I:I_r^t$. Then $I'R_Q=IR_Q$. Hence $Q\in\Sing(R/I')$. So we obtain $$\Sing(R/I)=V(J_1)\cup\Sing(R/I').$$ Note that $\height I' > r$ because $I'$ is the intersection of all primaty components of $I$ whose associated primes do not contain $I_r$.\par
By [T2, Lemma 1.1], $\height I_\alpha= r$. By Lemma 2.4, every generating set of  an ideal of $R$ specialize to a set of generators of the specialized ideal. Hence $A_\alpha$ is a Jacobian matrix of $I_\alpha$ and $(J_1)_\alpha$ is the ideal generated by the $r\times r$ minors of $A_\alpha$. Let $P_1,\ldots,P_s$ be the associated prime ideals of $I$ of height $r$. We have $(I_r)_\alpha=(P_1)_\alpha\cap\ldots\cap (P_s)_\alpha$ by Lemma 2.3. By [Kr2, Satz 14] and  [Se, Appendix, Theorem 7], each $(P_j)_\alpha$ is the intersection of prime ideals of height $r$. Hence these prime ideals are associated prime ideals of $I_\alpha$. On the other hand, since $I'=I:I_r^t$ for $t$ large enough, we have $I'_\alpha=\cup_{t\ge 1}(I:I_r^t)_\alpha=\cup_{t\ge 1}I_\alpha:(I_r)^t_\alpha$ by Lemma 2.5. Therefore,  $I'_\alpha$ is the intersection of all primary components of $I_\alpha$ which do not contain $(I_r)_\alpha$. But $\height I'_\alpha=\height I' > r$ by [T2, Lemma 1.1]. Thus, every associated prime ideal of $I_\alpha$ of height $r$ must contain $(I_r)_\alpha$. So $(I_r)_\alpha$ is the intersection of the associated prime ideals of $I_\alpha$ of height $r$. Similarly as above, we can show that
$$\Sing(R_\alpha/I_\alpha)=V((J_1)_\alpha)\cup\Sing(R_\alpha/I'_\alpha).$$\par
Now, using backward induction on $r$ we may assume that there is an ideal $J_2\supseteq I'$ such that $\Sing(R/I')=V(J_2)$ and $\Sing(R_\alpha/I'_\alpha)=V((J_2)_\alpha)$. Summing up we obtain\md
$\hspace{4cm} \Sing(R/I)=V(J_1)\cap V(J_2) = V(J_1+J_2),$\par $\hspace{3cm} \Sing(R_\alpha/I_\alpha)=V((J_1)_\alpha)\cap V((J_2)_\alpha)=V((J_1+J_2)_\alpha).$ \end{pf} 

Let $t \ge 0$ be a fixed integer. A  ring $A$ is said to satisfy  the Serre's condition $(S_t)$ if $\depth A_{\frak p} \ge \min\{\dim A_{\frak p},t\}$ for any prime ideal $\frak p$ or the Serre's condition $(R_t)$ if $A_{\frak p}$ is regular for any prime ideal $\frak p$ with $\height p\le t$. 

\begin{Theorem} Let $k$ be a field of characteristic zero. Let $\frak a$ be an ideal of $S$. Assume that $S/\frak a$ satifies one of the following properties:\sk
{\rm (i)} $(S_t)$,\sk {\rm (ii)} $(R_t)$,\sk {\rm (iii)} $S/\frak a$ is reduced, \sk {\rm (iv)} $S/\frak a$ is normal,\sk {\rm (v)} $S/\frak a$ is regular.\sk
\noindent Then $S_\alpha/{\frak a}_\alpha$ has the same property for almost all $\alpha.$\end{Theorem}

\begin{pf} Let $\frak b$ be an ideal of $S$ such that $\NCM(S/{\frak a})=V({\frak b})$ and $\NCM(S_\alpha/{\frak a}_\alpha)=V({\frak b}_\alpha)$ as in Lemma 4.3. Let ${\frak p} \supseteq \frak b$ be an arbitrary primde ideal of $S$. Then $\depth(S/{\frak a})_{\frak p}<\dim(S/{\frak a})_{\frak p}$. If $S/\frak a$ satisfies $(S_t)$, then $\depth (S/{\frak a})_{\frak p}\ge t$. Hence $\grade({\frak b}/{\frak a})\ge t$ [BH, Proposition 1.2.10]. By Corollary 3.4, $\grade({\frak b}_\alpha/{\frak a}_\alpha)\ge t$. Hence $S_\alpha/{\frak a}_\alpha$ satisfies $(S_t)$. Similarly, using Lemma 4.4 we can find an ideal $\frak c$ in $S$ such that $\Sing(S/{\frak a})=V({\frak c})$ and $\Sing(S_\alpha/{\frak a}_\alpha)=V({\frak c}_\alpha)$. If $S/\frak a$ satisfies $(R_t)$, then $\height{\frak c}/{\frak a}>t$. By Corollary 2.7, $\height{\frak c}_\alpha/{\frak a}_\alpha>t$. Hence $S_\alpha/{\frak a}_\alpha$ satisfies $(R_t)$.\par
The statement on the reduced property (normality) follows from the fact that a ring is reduced (normal) if and only if it satisfies $(S_1)$ and $(R_0)$ ($(S_2)$ and $(R_1)$). Finally, if $S/\frak a$ is regular, then $S/\frak a$ satisfies $(R_d)$, $d=\dim S/\frak a$. Hence $S_\alpha/{\frak a}_\alpha$ also satisfies $(R_d)$. Since $\dim S_\alpha/{\frak a}_\alpha=d$ by Theorem 2.6 (ii), this implies that $S_\alpha/{\frak a}_\alpha$ is regular.\end{pf}

\noindent{\it Remark.} One can apply the above results to study hypersurface sections. Indeed, that amounts to study rings of the form $k[X]/(I,f_\alpha)$, where $I$ is an ideal of $k[X]$ and $f_\alpha= \alpha_1f_1+\cdots+\alpha_mf_m$ with given $\alpha_1,\ldots,\alpha_m \in k$ and $f_1,\ldots,f_m \in k[X]$. These rings are specializations of the ring $k(u)[X]/(I,f)$, where $u=(u_1,\ldots,u_m)$ is a family of indeterminates and $f=u_1f_1+\cdots+u_nf_n$. Now we only need to find conditions on $I$ which allows a given property  to be transferred from $I$ to $(I,f)$ (see [H], [T1]) and to check whether this property is preserved by specializations. By this way we can reprove Bertini theorems from [Fl].\md
 
\section*{References}\md
 
\noindent [BH]  W. Bruns and J. Herzog, {\it Cohen-Macaulay rings,} Cambridge University Press, 1993.\par
\noindent [BE]  D.A. Buchsbaum and D. Eisenbud, {\it What makes a complex exact ?,} J. Algebra 25 (1973), 259-268.\par
\noindent [E]   D. Eisenbud, {\it Commutative algebra with a view toward algebraic geometry, } Springer-Verlag, 1995.\par
\noindent [Fa]  G. Faltings, {\it\"Uber die Annulatoren lokaler Kohomologiegruppen,} Archiv der Math. 30 (1978), 473-476.\par
\noindent [Fl] H. Flenner, {\it Die S\"atze von Bertini f\"ur lokale Ringe}, Math. Ann. 229 (1977), 97-111.\par
\noindent [H] M. Hochster, {\it Properties of noetherian rings stable under general grade reduction}, Archiv Math. 24 (1973), 393-396.\par
\noindent [Kr1]  W. Krull, {\it Parameterspezialisierung in Polynomringen,} Arch. Math.~1 (1948), 56-64.\par
\noindent [Kr2]  W. Krull, {\it Parameterspezialisierung in Polynomringen {\rm II}, Grundpolynom}, Arch. Math. 1 (1948), 129-137.\par
\noindent [Ku1]  W.E. Kuan, {\it A note on a generic hyperplane section of an algebraic variety,} Can. J. Math. 22 (1970), 1047-1054.\par
\noindent [Ku2]  W.E. Kuan, {\it On the hyperplane section through a rational point of of an algebraic variety,} Pacific. J. Math. 36 (1971), 393-405.\par
\noindent [Ku3]  W.E. Kuan, {\it Specialization of a generic hyperplane section through a rational point of an algebraic variety,} Ann. Math. Pura Appl. 94 (1972), 75-82.\par
\noindent [NT]  D.V. Nhi and N.V. Trung, {\it Specialization of modules,} to appear in Comm. Algebra.\par
\noindent [Sc]  P. Schenzel, {\it Einige Anwendungen der lokalen Dualit\"at und verallgemeinerte Cohen-Macaulay-Moduln,} Math. Nachr. 69 (1975), 227-242.\par  
\noindent [Se]  A. Seidenberg, {\it The hyperplane sections of normal varieties, } Trans. Amer. Math. Soc. 69 (1950), 375-386.\par
\noindent [T1] N.V. Trung, {\it\"Uber die \"Ubertragung der Ringeigenschaften zwischen $R$ und $R[u]/(F)$}, Math. Nachr. 92 (1979), 215-229.\par
\noindent [T2]  N.V. Trung, {\it Spezialisierungen  allgemeiner Hyperfl\"achenschnitte und Anwendungen}, in: Seminar D.Eisenbud/B.Singh/W.Vogel, Vol.~1, Teubner-Texte zur Mathematik, Band 29 (1980), 4-43.\par
\noindent [T3]  N.V. Trung, {\it Toward a theory of generalized Cohen-Macaulay modules,} Nagoya Math. J. 102 (1986), 4-43.
\end{document}